\documentclass[a4paper,12pt]{amsart}
\usepackage{multicol,ifthen}
\usepackage{amssymb,stmaryrd}
\usepackage{harvard}
\usepackage{graphicx}

\newcommand{\R}{\mathbb R}

\newcommand{\C}{\mathbb C}

\newcommand{\beq}{\begin{equation}}
\newcommand{\eeq}{\end{equation}}
\newcommand{\beqarr}{\begin{eqnarray}}
\newcommand{\eeqarr}{\end{eqnarray}}
\newcommand{\beqa}{\begin{eqnarray*}}
\newcommand{\eeqa}{\end{eqnarray*}}
\begin{document}
\hyphenation{Rie-man-nian}

\thispagestyle{empty}
\renewcommand{\thefootnote}{\fnsymbol{footnote}}
\title[Weyl geometry  in late 20th century physics]{Weyl geometry in late 20th century physics}
\author[E. Scholz]{Erhard Scholz}

\date{23. 09. 2009  }
\renewcommand{\thefootnote}{\arabic{footnote}}

\begin{abstract}
Weyl's original scale geometry of 1918 (``purely infinitesimal geometry'') was withdrawn from physical theory in the early 1920s. It had a comeback in the last third of the 20th century in different contexts: scalar tensor theories of gravity, foundations of physics (gravity, quantum mechanics), elementary particle physics, and cosmology.  Here we survey the last two segments. It seems that Weyl geometry continues to have an open research potential for the foundations of  physics after the turn of the century.
\end{abstract}

\maketitle

\section{Introduction}

Roughly at the time when his famous book  {\em Raum $\cdot$ Zeit $\cdot$ Materie} (RZM) went into print, Hermann Weyl generalized 
  Riemannian geometry by introducing scale freedom of the underlying metric, in order to bring a  more basic ``purely  infinitesimal'' point of view to bear \cite{Weyl:InfGeo,Weyl:GuE}. How Weyl extended his  idea of scale gauge to a unified theory of the electromagnetic and gravitational fields, how this proposal was  received among physicists,  how it was given up   -- in its original form  --  by the inventor already two years later, and how it was transformed into  the now generally accepted  $U(1)$-gauge theory of the electromagnetic field, has been extensively 
studied.\footnote{\cite{Vizgin:UFT,Straumann:Einstein_Weyl,Skuli:Diss,Goenner:UFT,ORaif/Straumann,%
ORaifeartaigh:Dawning,Scholz:DMV,Scholz:Weyl_PoS,Scholz:Fock/Weyl}.} 
 Many times Weyl's original  scale gauge geometry was  proclaimed   dead, physically misleading or, at least, useless as a physical concept. But it  had  surprising come-backs in various research programs of physics. It seems well alive  at the turn to the new century.

Weylian  geometry was taken up explicitly or  
 half-knowingly  in different  research fields of theoretical physics during the second half of the 20th century  (very rough time schedule): 
\begin{itemize}
\item 1950/60s:  Jordan-Brans-Dicke theory
 \item 1970s: a double  retake of Weyl geometry by    Dirac and Utiyamah
\item   1970/80s: Ehlers-Pirani-Schild and successor studies
\item 1980s: geometrization of  (de Broglie Bohm) quantum potential 
\item 1980/90s: scale invariance and the Higgs mechanism
\item 1990/2000s: scale covariance in recent cosmology
\end{itemize}
 All these topics are  worth of closer historical studies. Here we concentrate on the last two topics. The first four have to be left to a more extensive study.

With the rise of the standard model of elementary particles (SMEP) during the 1970s a new context for the discussion of fundamental questions in  general relativity formed.\footnote{It is complemented by the standard model of cosmology, SMC. Both, SMC and SMEP, developed a peculiar symbiosis since the 1970s \cite{Kaiser:Colliding}. Strictly speaking, {\em the} standard model without further specifications consists of the two closely related complementary parts SMEP and SMC. }  
That led to an   input of new ideas  into gravity. Two subjects played a crucial role for our topic:  {\em scale} or {\em   conformal invariance} of the known interactions of high energy physics (with exception of gravity) and  the intriguing idea of {\em symmetry reduction}  imported  from solid state physics to the electroweak sector  of the standard model. The latter is usually understood as  symmetry breaking due to some  dynamical process (Nambu, Goldstone, Englert, Higgs, Kibble e.a.). The increasingly successful standard model worked  with conformal invariant  interaction fields, mathematically spoken connections with values in the Lie algebras of ``internal''  symmetry groups (i.e., unrelated to the spacetime), $SU(2)\times U(1)_Y$ for the electroweak (ew) fields, $SU(3)$ for the chromodynamic field modelling  strong interactions, and  $U(1)_{\text{em}}$ for the electromagnetic (em) field,  inherited from the 1920s.   In the SMEP electromagnetism appears as a  residual phenomenon, after breaking  the  isospin  $SU(2)$   symmetry of  the { ew} group  to the isotropy group $U(1)_{\text{em}}$ of a hypothetical vacuum state. The latter is  usually characterized by  a {\em Higgs field} {\bf $ \Phi$}, a  ``scalar'' field  (i.e. not transforming under spacetime coordinate changes) with values in an isospin  $\frac{1}{2}$ representation of the weak $SU(2)$ group. If    $ \Phi$ characterizes  dynamical  symmetry breaking, it should have a  massive quantum state, the Higgs boson  \cite{Higgs:1964,Weinberg:1967}. The whole procedure became  known under the  name  
``Higgs mechanism''.\footnote{Sometimes called in more length and greater historical justness ``Englert-Brout-Higgs-Guralnik-Hagen-Kibble'' mechanism.}

Three interrelated questions arose naturally if one wanted to bring gravity closer to the physics of the standard model: 
\begin{itemize}
\item[(i) ] Is it possible to bring conformal, or at least scale covariant generalizations of classical (Einsteinian) relativity 
 into a coherent common  frame with the standard model SMEP?\footnote{Such an attempt seemed  to be  supported experimentally by the phenomenon of (Bjorken) scaling in deep inelastic
 electron-proton scattering experiments. The latter indicated, at first glance, an active scaling symmetry of mass/energy in high energy physics; but it turned out to hold only approximatively and of restricted range.}
\item[(ii)] Is  it  possible to embed classical relativity in a quantized theory of gravity? 
\item[(iii)] Or just the other way round,  can ``gravity do something like  the Higgs''?\footnote{Formulation due to 
 \cite{Pawlowski:1990}.} That would be the case if the mass acquirement of electroweak bosons could  be understood by a Brans-Dicke like extension of gravitational structures.
\end{itemize}
These questions were posed and attacked with differing degrees of success since the 1970s to the present. Some of these contributions, mostly referring to questions (i) and (iii),  were closely related to Weylian scale geometry or even openly formulated in this framework. The literature on these questions is immense. Obviously we can only scratch on the surface of it in our survey, with strong selection according to the criterion  given by the title of this paper. So we exclude discussion of topic (ii), although it was historically closer related to the other ones than it appears here (section \ref{section Higgs}).

 In the last three decades of the 20th century a dense cooperation between  particle physics, astrophysics and cosmology was formed. The emergence of this   intellecual and disciplinary symbiosis had many causes; some of them are discussed in \cite{Kaiser:Colliding} and  by C. Smeenk (this volume). Both papers share a common interest in inflationary theories of the very early universe. But in the background of this reorganization  more empirically driven changes, like  the accumulating evidence for   ``dark matter'' by astronomical observations  in the 1970s, were surely of great importance \cite{Rubin:dark_matter,Trimble:dark_matter}. That had again strong  theoretical   repercussions. In the course of the 1980/90s it forced astronomers and astrophysicists to assume a large amount of  non-visible, 
non-baryonic matter with rather peculiar properties.   In the late 1990s  increasing and different evidence spoke strongly in favour of a   non-vanishing cosmological constant  $\Lambda $. It was now  interpreted as a  ``dark energy'' contribution to the dynamics of the universe \cite{Earman:Lambda}.  

 The second part of the 1990s  led to a relatively coherent picture of the {\em standard model of cosmology} SMC with a  precise  specification of the values of the energy densities $\Omega _m,\Omega _{\Lambda }$ of (mostly ``dark'') matter and of ``dark'' energy as the central parameters of the model. This  specification depended, of course, on the choice of the Friedman-Lemaitre spacetimes as theoretical reference frame.  $\Omega _m$ and $\Omega _{\Lambda }$ together determine the adaptable parameters of this model class (with cosmological constant). The result was the now favoured  
$\Lambda $CDM model.\footnote{CDM stands for for cold dark matter and $\Lambda $ for  a non-vanishing cosmological constant.}
In this sense, the geometry of the physical universe, at least its empirically accessible part, seems to be well determined, in distinction to the quantitative underdetermination  of many of the earlier cosmological world pictures of extra-modern or early modern cultures \cite{Kragh:Cosmos}. But the new questions related to ``dark matter'' and ``dark energy'' also induced attempts for  widening  the frame of classical GRT. Scale covariant scalar fields in the framework of conformal geometry, Weyl geometry, or Jordan-Brans-Dicke theory formed an important cluster of such alternative attempts. We shall have a  look at them in section (\ref{section cosmology}). But before we enter this discussion, or pose the question of the role of scale covariance in particle physics, we give a short review of the central features of Weyl geometry, and its relation to Brans-Dicke theory, from a systematic point of view. Readers with a background in these topics might like to skip the next section and pass directly to section (\ref{section Higgs}). The paper is concluded by a short evaluation of our survey (section \ref{section discussion}).

\section{Preliminaries on Weyl geometry and  JBD theory \label{section preliminaries}}
\subsection*{Weylian metric, Weyl structure}
A  Weylian metric on a differentiable manifold $M$ (in the following mostly $dim \; M = 4$) can be given by pairs $(g, \varphi)$ of a non-degenerate symmetric differential two form $g$, here of Lorentzian signature  $(3,1)= (-,+,+,+,)$, and a differential 1-form 
 $\varphi$.
 The {\em Weylian metric}  consists of the  equivalence class of such pairs, with $(\tilde{g} , \tilde{\varphi} ) \sim (g,   \varphi) $  iff
\beq  (i) \quad  \tilde{g} =  \Omega ^2 g \; , \quad \quad  (ii) \quad \tilde{\varphi} = \varphi - d \log \Omega \, \label{gauge transformation}\eeq 
for  a strictly positive real function  $\Omega > 0$  on $M$. Chosing a representative means to {\em gauge} the Weylian metric;   $g$ is then  the {\em Riemannian component}  and $\varphi$ the {\em scale connection} of the gauge. A change of representative (\ref{gauge transformation}) is  called a {\em Weyl} or {\em scale transformation}; it consists of a  conformal rescaling $(i)$ and a {\em scale gauge transformation} $(ii)$. 
A manifold with a Weylian metric $(M, [g,\varphi ])$ will be called a  {\em Weylian manifold}. 
For more detailed introductions to Weyl geometry in the theoretical physics literature see \cite{Weyl:RZM,Bergmann:Relativity,Dirac:1973}, for mathematical introductions \cite{Folland:WeylMfs,Higa:1993}. 

In the recent mathematical literature
 a {\em Weyl structure}  on a differentiable manifold $M$ is specified by  a pair $(c,\nabla)$  consisting of  a conformal structure $c=[g]$ and an affine, i.e.  torsion free, connection $\Gamma$, respectively its covariant derivative $\nabla$. The latter is constrained  by the property that for any $g\in c$ there is a differential 1-form $\varphi_g$ such that 
\beq \nabla g + 2 \varphi_g \otimes g  = 0 \, , \label{compatibility condition Weyl structure} \eeq
 \cite{Calderbank/Pedersen,Gauduchon:1995,Higa:1993,Ornea:2001}. We shall call this {\em weak compatibility} of the affine connection with the metric.\footnote{Physicists usually prefer to speak of a ``semimetric connection''' \cite{Hayashi/Kugo:1977} or even of a ``nonmetricity'' of the connection \cite{Hehl_ea:Report1995} etc.}
 One could also formulate the compatibility by
 \beq  \Gamma -  _g\hspace{-0.2em}\Gamma   = 1 \otimes \varphi_g + \varphi_g \otimes 1 - g \otimes \varphi_g^{\ast}\, ,  \label{2nd compatibility Weyl structure}  \eeq
 where $1$ denotes the identiy in $Hom(V,V)$ for every $V=T_xM$,   $\varphi_g^{\ast}$ is the dual of $\varphi_g$ with respect to $g$, and $ _g\hspace{-0.1em}\Gamma $ is the Levi-Civita conection of $g$. Written in coordinates that means 
\beq  \label{Christoffel}  \Gamma^{\mu }_{\nu \lambda } =   {}_g\Gamma^\mu _{\nu \lambda } + \delta ^{\mu }_{\nu } \varphi _{\lambda } +
\delta ^{\mu }_{\lambda } \varphi _{\nu } - g_{\nu \lambda } \varphi^{\mu } \, , \eeq 
if  $ {}_g\Gamma^\mu _{\nu \lambda }$ denote the coefficients of the affine connection with respect to the Riemannian component $g$ only.

This is just another way to specify the structure of a Weylian manifold, because  $[(g,\varphi )]$  is compatible with exactly one affine connection.  (\ref{scale covariant derivative of g}) is the condition that the scale covariant derivative of $g$ vanishes in every gauge (see below).
The Weyl structure is called {\em closed}, respectively {\em exact}, iff the differential 1-form $\varphi_g$ is so (for any $g$). In agreement with large parts of the physics literature on Weyl geometry, we shall  use the terminology {\em integrable} in the sense of closed, i.e., in a local sense.

In some part of the physics literature a {\em change of scale}  like in (\ref{gauge transformation}\,$(i)$) is considered without explicitly mentioning  the  accompanying {\em gauge transformation} $(ii)$. Then a scale transformation is identified with  a conformal transformation of the metric. That may be misleading but need not, if the second part of (\ref{gauge transformation}) is respected  indirectly. 
In any case we have to distinguish between a {\em strictly conformal} point of view and a {\em Weyl geometric} one. In the first case we deal with $c=[g]$ only, in the second case we refer to the whole Weyl metric $[(g,\varphi)]$, respectively Weyl structure $(c,\nabla)$.\footnote{Both approaches work with the ``localized'' (physicists' language)   {\em scale extended Poincar\'e group} $\mathcal{W}=\R^4 \ltimes SO^+(3,1)\times \R^+$ as gauge automorphisms. The {\em transition} from a strictly conformal  approach to a  Weyl geometric one {\em has nothing to do with a group reduction} (or even with ``breaking'' of some symmetry); it rather  consists of an enrichment of the structure  while upholding the automorphism group. } 

\subsection*{Covariant derivative(s), curvature, Weyl fields}
The  covariant derivative with respect to $\Gamma$ will be denoted  (like above)  by $\nabla  $. The covariant derivative with respect to the Riemannian component of the metric only will be indicated by ${}_g\hspace{-0.1em}\nabla $.  $\nabla $ is  an invariant operation for  vector and tensor fields on $M$, which are themselves  invariant under gauge transformations. The same can be said for {\em geodesics} $\gamma_W$ of Weylian geometry, defined by $\nabla $,  and for the  {\em Riemann curvature}  tensor $Riem = (R^{\alpha} _{\beta \gamma \delta})$ and its contraction, the Ricci tensor $Ric= (R_{\mu \nu})$. The contraction is defined with respect to the 2nd and 3rd component
\beq R_{\mu \nu} := R^{\alpha} _{\mu \alpha \nu}
\eeq .

Functions or (vector, tensor, spinor \ldots) fields  $F$ on $M$, which transform under gauge transformations like
\beq   F \longmapsto \tilde{F} = \Omega ^k F  \; . \label{Weyl fields} \eeq
will be called {\em Weyl functions} or {\em Weyl fields} on $M$ of  (scale or Weyl) {\em weight} 
$  w(f):= k$. Examples are: $w(g_{\mu \nu })=2$, $w(g^{\mu \nu })=-2$ etc.  As the curvature tensor $Riem$  of the Weylian metric and  the Ricci curvature $Ric$ are scale invariant,      scalar curvature 
\[   {R} : = g^{\mu \nu } R\, _{\mu \nu } \] 
is  of weight $ w({R}) =-2$. For the sake of historical precision it has to be noted that Weyl himself  considered   $g$  to be of weight $\overline{w} (g)  = 1$. Accordingly Weyl's original weights, and those of a   considerable part of the literature, are half of ours,  $\overline{w} = \frac{1}{2}w$. Moreover, in most of the physics literature the sign convention for the scale connection is different; both together means that a differential form $\kappa = - 2 \varphi $ is used in the description of  Weyl geometry.\footnote{Reasons for our  conventions: Our  sign choice of the scale connection implies positive exponent of the scale transfer function (\ref{scale transfer}). Our   weight convention is such that the   length (norm) of  vectors has  weight  $1$.    } 

The covariant derivative $\nabla$ of Weyl fields $F$ of weight $w(F) \neq 0$   does not lead to a  scale covariant quantity. This is a  deficiency of the geometric  structure considered so far, if one works  in a field theoretic  context. It can be repaired by introducing a {\em scale covariant derivative} $D$ of Weyl  fields in addition to the scale invariant  $\nabla $:
\beq D F := \nabla F  + w(F)
\varphi \otimes F \,.    \label{scale covariant derivative}  \eeq
A scale covariant vector field $F^{\nu }$, e.g., has the scale covariant derivative 
\[  D_{\mu } F^{\nu } := \partial _{\mu } F^{\nu } + \Gamma^{\nu }_{\mu \lambda } F ^{\lambda} + w(F)\, \varphi_{\mu } F^{\nu }  \, , \]
with the abbreviation
$\partial_{\mu } := \frac{ \partial }{\partial x^{\mu } }$ etc. The compatibility condition in the definition of a Weyl structure (\ref{compatibility condition Weyl structure}) can now be written as
\beq  D g = 0  \label{scale covariant derivative of g}  \, . \eeq

\subsection*{Relation to Jordan-Brans-Dicke theory}
 Jordan-Brans-Dicke   (JBD)  theory assumes  a scalar field $\chi$ of scale weight $w(\chi)=-1$, coupled to gravity (a pseudo-Riemannian metric $g$) by a Lagrangian of the following type
\beq \mathcal{L_{JBD}}(\chi, g) = (  \chi {R} -  \frac{\omega}{\chi }\partial ^{\mu} \chi\, \partial _{\mu} \chi  ) \sqrt{|det\, g|} \, , \label{Lagrangian JBD} \eeq
with a free parameter $\omega$ and scalar curvature $R$. It considers conformal transformations of metric and fields, while fixing the  Levi-Civita connection $\nabla$ of the metric $g$ underlying   (\ref{Lagrangian JBD}). Such a  conformal rescaling is called a change of {\em frame}. The ``original'' one (defining the affine connection as the Levi-Civita connection of the Riemannian metric) like in (\ref{Lagrangian JBD}) is called  {\em Jordan frame}. The  one in which the scalar field  (und thus  the coefficient of the Einstein-Hilbert term, the gravitational coupling coefficient)  is scaled to  a constant is called {\em Einstein frame}. 

A conformal class of a metric $[g]$ and specificiation of an affine connection like in JBD theory characterizes an integrable Weyl structure.  We should thus be aware that  JBD theory carries the basic features of a Weyl geometric structure, even though most of the workers in the field do not look at it from this point of view.  In this sense I  consider JBD theory  as a research field in which Weyl geometry stood in the background ``half-knowingly''.\footnote{One need not know Weyl geometry, in order to work in the framework of such a naturally given  structure;  just like Moli\'ere's M. Jourdain did not know that he had spoken  prose for forty years, before he was told so by a philosopher.}   
Jordan and Einstein frame are nothing but  Riemann gauge, respectively scalar field gauge, in the language just introduced for  integrable Weyl geometry.

More recent literature, like the excellent monographs \cite{Fujii/Maeda,Faraoni:2004}, often  prefers a  slightly different form of the scalar field and the Lagrangian, $\phi = \sqrt{2 \xi^{-1}\chi} $ (scale weight  $w(\phi)= -1$),  $\xi=  \frac{\epsilon }{4 \omega }$. Then the Lagrangian  acquires the form 
\beq  \mathcal{L}_{BD}= \left(\frac{1}{2} \xi  \phi^2 {R} - \frac{1}{2}\epsilon  \partial ^{  \mu }  \phi \, \partial _{\mu}\phi  + L_{mat} \right) \sqrt{|det\, g|}     \, ,\label{Lagrangian JBD Fujii/Maeda} \eeq 
where  ${sig}\, g = (3,1) \cong  (-+++)$ and   $\epsilon = \pm 1$ or  $0$ \cite[5]{Fujii/Maeda}.\footnote{ According to Fujii/Maeda  $\epsilon =1 $ corresponds to a ``normal field having a positive energy, in other words, not a ghost''.  $\epsilon =-1$  may look at first unacceptable because it ``seems to indicate negative energy'', but ``this need not be an immediate difficulty owing to the presence of the nonminimal coupling.'' '(ibid.)} 
 \citeasnoun{Penrose:1965} showed that  $ \mathcal{L}_{BD}$ is conformal invariant for $\xi = \frac{n-2}{4(n-1)}$ ($n$ spacetime dimension). 

Moreover in the recent literature strong arguments have accumulated to prefer Einstein gauge over Jordan gauge \cite{Faraoni_ea:Transformations}. An obvious argument comes from the constraints of the coefficient $\omega$ arising from high precision gravity observations in the solar system, if Jordan frame is considered to be ``physical''.\footnote{See the contribution by C. Will, this volume.}

\section{Scale covariance in   particle physics \label{section Higgs}}
\subsection*{Englert's conformal approach}
Fran\c{cois} Englert and  coworkers   studied 
 conformal gravity as part of the quantum field program \cite{Englert_ea:1975}. In a common paper written with the astrophysicist  Edgar Gunzig and others, the authors established an explicit link to JBD theory (not  to  Weyl geometry). They 
started from a ``dimensionless'', i.e. scale invariant,  Lagrangian for gravitation with a square curvature term of an affine connection  $\Gamma$ {\em not} bound to the metric,  $\mathcal{L}_{\text{grav}} = R^2 \, \sqrt{|det\, g|} $  in addition to a Lagrangian matter term  \cite{Englert_ea:1975}. In consequence, the authors  varied with respect to the metric $g$ and the connection $\Gamma$ {\em independently}. 

Further compatibility considerations made the connection weakly metric compatible, in the sense of our equ.  (\ref{compatibility condition Weyl structure}), even with an integrable scale connection \cite[equ.(7)]{Englert_ea:1975}. In this way, the approach worked in a Weyl structure, but the authors did not care about it. They rather tried to be as ``conformal'' as possible.

  In an attempted  ``classical phenomenological description'' they  characterized a pseudo-Riemanian  Lagrangian  of a scalar field coupled to gravity like in our equ. (\ref{Lagrangian JBD Fujii/Maeda}),  with the necessary specification $\xi = \frac{1}{6}$ in order to achieve conformal symmetry.  The  scalar field was called  ``dilaton'' and  considered  as  a
 ``Nambu-Goldstone boson''  of a  ``dynamical symmetry breakdown'' of the scale symmetry, but {\em without} a massive ``scalar meson'' \cite[75]{Englert_ea:1975}. The terms corresponding  to the Weylian scale connection (re-reading their paper in the light of Weyl geometry) were not considered as a physical field, but as a mathematical artefact of the analysis. 

In one of the following papers Englert, now   with other coauthors,  studied the perturbative behaviour of  conformal gravity  ($\xi = \frac{n-2}{4(n-1)}$) coupled to massless fermions and photons in $n \geq 4$ dimensions.\footnote{The motivation to consider $n \geq 4$ was dimensional regularization.} They came to the conclusion that anomalies arising in the calculations for non-conformal actions disappeared at the tree and 1-loop levels in their approach. They took this as an indicator that  gravitation might perhaps arise in a ``natural way from spontaneous breakdown of conformal invariance''  \cite[426]{Englert_ea:1975}. 

\subsection*{Smolin introduces  Weyl geometry}
Englert's e. a. paper was one of the early steps into the direction (i) of our introduction. 
Other authors followed and extended this view, some of them explicitly in a Weyl geometric setting, others clothed in the language of conformal geometry. The first strategy was  chosen by  Lee Smolin in his paper \cite{Smolin:1979}. In section 2 of the paper he gave an explicit and clear introduction to Weyl geometry.\footnote{In his bibliography he went back directly to \cite{Weyl:RZMEnglish} and \cite{Weyl:GuE}; he  did not   quote any of the later literature on Weyl geometry. }
The  ``conformally metric gravitation'', as he called it, was built  upon a matter-free Lagrangian with Weyl geometric curvature terms $R, \, Ric = (R_{\mu \nu })$,  $f=(f_{\mu \upsilon })$ for scale curvature alone, and scale covariant Weylian derivatives $D$ (in slight adaptation of notation): 
\beqarr  |det\, g|^{-\frac{1}{2}}  \mathcal{L}_{\text{grav}} &=&  -\frac{1}{2}c \,  \phi^2 R + \; [  - {e_1} R^{\mu \nu }R_{\mu \nu } -  {e_2}  R^2 ] \label{Smolin's Lagrangian}\\
& & + \quad \frac{1}{2}D^{\mu}\phi \,D_{\mu}\phi -\frac{1}{4g^2}f_{\mu \nu } f^{\mu \nu } - \lambda \phi^4 \nonumber
 \eeqarr 
 with coupling coefficients $c, e_1, e_2, g, \lambda $.\footnote{Signs  have to be taken with caution. They may depend on conventions for defining the Riemann curvature, the Ricci contraction, and the signature. Smolin, e.g., used a different sign convention for $Riem$  to the one used in this survey. Signs given here are adapted to ${signature}\, g = (3,1)$, Riemann tensor of mathematical textbooks, and Ricci contraction like in section \ref{section preliminaries}. \label{fn signs} }
For  coefficients 
of the quadratic curvature terms (in square brackets) with  $e_2 = - \frac{1}{3}e_1 $, the latter  was variationally equivalent (equal up to divergence) to the  squared conformal curvature $C^2 = C_{\mu \nu \kappa \lambda } C^{\mu \nu \kappa \lambda }.$\footnote{General knowledge, made explicit, e.g. by \cite{Hehl_ea:quad_curv}.}

Smolin introduced the scalar field $\phi$ not only by formal reasons (``to write a conformally invariant Lagrangian with the required properties''), but with similar physical  interpretations as Englert e.a.,\footnote{\cite{Englert_ea:1975} was not quoted by Smolin.}
namely ``as an order parameter to indicate the spontaneous breaking of the conformal invariance'' \cite[260]{Smolin:1979}. His Lagrangian used  a modified adaptation from JBD theory, ``with some additional couplings'' between scale connection  $\varphi$ and scalar field $\phi$. But Smolin  emphasized that ``these additional couplings go against the spirit of
 Brans-Dicke theory''  as they introduced a non-vanishing divergence of the non-gravitational fields.

For low energy considerations  Smolin dropped the square curvature term (in square brackets,  (\ref{Smolin's Lagrangian})),  added an ``effective'' potential term of the scalar field $V_{\text{eff}}(\phi )$ and derived the equations of motion by varying with respect to  $g, \phi, \varphi$. Results were  Einstein equation,  scalar field equation, and Yang-Mills equation for the scale connection.

 Smolin's  Lagrangian contained  terms in the scale connection:\footnote{Smolin's complete Lagrangian was  \[  |det\, g|^{-\frac{1}{2}}  \mathcal{L}^{grav} =  \frac{1}{2}c \, F^2\,  _g\hspace{-1pt}R -\frac{1}{4g^2}f_{\mu \nu } f^{\mu \nu } + \frac{1}{8}(1+6c)F^2 \varphi_{\mu} \varphi^{\mu } - V_{\text{eff}}(F) \; . \] 
}
\beq  -\frac{1}{4g^2}f_{\mu \nu } f^{\mu \nu } + \frac{1}{8}(1+6c)F^2 \varphi_{\mu} \varphi^{\mu }
\eeq 
That looked like a mass term for  $\varphi$ considered as  potential of the scale curvature field $f_{\mu \nu }$, called ``Weyl field'' by Smolin. By comparison with the Lagrangian of the Proca equation in electromagnetic theory,  Smolin   concluded that the ``Weyl field'' has  mass close to the Planck scale, given by
\beq  M^2_{\varphi}  =  \frac{1}{4}(1+6c)F^2 \, . \eeq 
He commented that in his Weyl geometric gravitation theory  ``general relativity couples to a massive vector field'' $\varphi$. The scalar field $\phi$, however, ``may be absorbed into the scalar parts'' of $g_{\mu \nu }$ and $\varphi_{\mu}$ by a change of variables and remains massless \cite[263]{Smolin:1979}. In this way, 
Smolin brought Weyl geometic gravity closer to the field theoretic frame of particle physics. He did not discuss mass and interaction fields of the SMEP. Morover, the huge mass of the ``Weyl field'' must have appeared quite irritating.

\subsection*{Interlude}
At the time Smolin's paper appeared,  the program of  so-called {\em induced gravity}, entered an  active phase. Its central goal was to derive the action of conventional or modified Einstein gravity from an extended scheme of standard model type quantization. Among the authors involved in this program Stephen Adler and Anthony Zee stick out.   We cannot go  into this story here.\footnote{ For a survey of the status of investigations in 1981 see  \cite{Adler:Report_1982}; but note in particular 
\cite{Zee:1982a,Zee:1983}. The topic of ``origin of spontaneous symmetry breaking'' by radiative correction was much older, see e.g. \cite{Coleman/Weinberg:1973}.
In fact, Zee's first publication  on the subject preceded Smolin's. \cite{Zee:1979} was submitted in December 1978 and  published in February 1979; \cite{Smolin:1979} was submitted in June 1979.}

 Smolin's view that already the  structure of Weyl geometry might be well suited  to bring classical gravity into a coherent frame with standard model physics did not find much direct response. But it was ``rediscovered''  at least twice (plus  an independently developed conformal version). In 1987/88   Hung Cheng at the MIT   and a  decade later  Wolfgang Drechsler and  Hanno Tann, both at Munich,   arrived  basically at similar insights. both with an explicit extension to standard model (SMEP) fields  \cite{Cheng:Vector_meson,Drechsler/Tann,Drechsler:Higgs}. Simultaneously to  Cheng, the core of the  idea was  once more discovered  by  Mosh\'e Flato (Dijon)  and Ryszard R\c{a}cka  (during that time at Trieste), although  they formulated it in a strictly conformal framework without  Weyl structure  \cite{Flato/Raczka}. 
Neither  Cheng, nor Flato/R\c{a}cka or Drechsler/Tann seem to have known Smolin's proposal (at least Smolin is not cited by them); even less did they refer to the papers of each other.\footnote{Flato/R\c{a}cka's paper appeared as a preprint of the {\em Scuola Internazionale Superiore di Studi Avanzati}, Trieste, in 1987; the  paper itself  was submitted in December 1987 to {\em Physics Letters B} and published in July 1988.   Cheng's paper was submitted in February 1988,  published in November. Only a  decade later, in March 2009, Drechsler and Tann got acquainted with the other two papers. This indicates  that the Weyl geometric approach in field theory has not yet acquired the coherence of a research program with a stable  subcommunity.  }
All three approaches had their own achievements. Here we can give only give a short presentation of the main points of the work directly related to Weyl geometry.

\subsection*{Hung Cheng and his  ``vector meson''}
Hung Cheng started out from a Weyl geometric background, apparently inherited from the papers of  Japanese authors around Utiyama. The latter had taken up Weyl geometry in the early 1970s in a way not too different from Smolin's later approach.\footnote{ \cite{Utiyama:1973,Utiyama:1975a,Utiyama:1975b,Hayashi/Kugo:1979} }
Hung Cheng   extended Utiyama's theory  explicitly to the electroweak sector of the SMEP. The scalar field  $\Phi $ of weight $-1$  (without a separate  potential) was supposed to have values in an  isospin $\frac{1}{2}$  representation.\footnote{In the sequel the isospin extended scalar field will  be denoted by $\Phi $.} Otherwise it coupled to  Weyl  geometric curvature  $R$ as known. 
\beqarr  \mathcal{L}_{R} &=&    \varepsilon \frac{1}{2} \beta \,  \Phi^{\ast} \Phi R  \, |det\, g|^{\frac{1}{2}}  \label{Hilbert term Weyl geometric} \\
 \mathcal{L}_{\Phi} &=&  \quad \frac{1}{2}\tilde{D}^{\mu} \Phi^{\ast}  \tilde{D}_{\mu}\Phi   \, |det\, g|^{\frac{1}{2}}  \; ,\label{phi term Weyl geometric}
\eeqarr
with $\varepsilon = 1$.\footnote{Drechsler and Tann would later find reasons to set $\varepsilon =-1$ (energy of the scalar field positive). Hung Cheng's curvature convention was not made explicit; so there remains a sign ambiguity.}
The  scale covariant derivatives were extended to a  ``localized" {\em ew} group $SU(2)\times U(1)$. With the usual denotation of  the standard model,   $W^j_{\mu} $ for the field components of the $su(2)$ part (with respect to the Pauli matrices $\sigma _j\; (j=0,1,2)$) and  $B_{\mu}$ for $u(1)_Y \cong \R$ and coupling coefficients $g, g'$ they read\footnote{Cheng added another coupling coefficient for the scale connection, which is here suppressed. }
\beq  \tilde{D}_{\mu}\Phi    = (\partial _{\mu} - \varphi_{\mu} + \frac{1}{2} i g W^j_{\mu} \sigma_j + \frac{1}{2} g'B_{\mu})\Phi \, . \label{ew covariant derivative}
\eeq 
Cheng added Yang-Mills interaction Lagrangians for {\em ew} interaction fields $F$ and $G$ of the potentials $W$ (values in $su_2$), respectively $B$ (values in  $u(1)_Y$),   and added a scalar curvature term in  $f=(f_{\mu \nu}) = d\varphi$
\beq   \mathcal{L}_{\text{YM}}= -\frac{1}{4}\left( f_{\mu \nu}  f^{\mu \nu}  + F_{\mu \nu}  F^{\mu \nu}  + G_{\mu \nu}  G^{\mu \nu}  \right) \, |det\, g|^{\frac{1}{2}} \, .  \label{Cheng's ew interaction Lagrangian} \eeq
Finally he introduced spin $\frac{1}{2}$ fermion fields $\psi$ with the  weight convention $w(\psi)=-\frac{3}{2}$, and a Lagrangian  $\mathcal{L}_{\psi}$ similar to the one formulated later by Drechsler,  discussed below   (\ref{Lagrangian Dirac field}).\footnote{The second term in (\ref{Lagrangian Dirac field}) is missing in Cheng's publication. That is probably not intended, but a misprint. Moreover he did not discuss scale weights for Dirac matrices in the tetrad approach.} 

Cheng called  the scale connection, resp. its curvature, {\em Weyl's meson} field. Referring to Hayashi's e.a. observation  that  the scale connection does not influence the equation of motion of the spinor fields, he  concluded:\footnote{Remember that the $\varphi$ terms of scale covariant derivatives in the Lagrangian of spinor fields cancel.}
\begin{quote}
\ldots  Weyl's vector meson does not interact with leptons or quarks. Neither does it interact with other vector mesons. The only interaction the Weyl's meson has is that with the graviton. \cite[2183]{Cheng:Vector_meson}
\end{quote}

 Because of the tremendous mass of ``Weyl's vector meson'' Cheng conjectured that even such a minute coupling might  be of some cosmological import. More precisely, he   wondered, ``whether Weyl's meson may account for at least part of the dark matter of the universe'' (ibid.). Similar conjectures were stated once and again over the next decades,  if theoretical entities were encountered which might represent massive particles without  experimental evidence. Weyl geometric field theory was not spared this experience.  

\subsection*{Can gravity do what the Higgs does?  }
In the same year in which Hung Cheng's paper appeared, Mosh\'e Flato and Ryszard  R\c{a}czka sket\-ched  an approach  in which they  put gravity into a quantum physical perspective.\footnote{More than a decade earlier Flato had worked out a covariant (``curved space'') generalization of the Wightman axioms  \cite{Flato/Simon}, obviously different from the one discussed by R. Wald in this volume, with another coauthor. }
Although it would be interesting to put this paper in perspective of point (ii) in our introduction, we cannot do it here. In our contect, this paper matters because it introduced a  scale covariant Brans-Dicke like field in an isospin representation similar to Hung Cheng's, but in a strictly conformal framework  \cite{Flato/Raczka}. 

Six years later, R. R\c{a}czka took up the thread again, now in  cooperation  with Marek 
Paw{\l}owski.   In the meantime Paw{\l}owski 
 had joined the research program  by a  paper  in which he addressed the question whether perhaps gravity ``can do what the Higgs does'' \cite{Pawlowski:1990}. In  a couple of  preprints
 \cite {Pawlowski/Raczka:1994a,Pawlowski/Raczka:1995_0,%
Pawlowski/Raczka:1995_I,Pawlowski/Raczka:1995_II} and two   refereed papers 
\cite{Pawlowski/Raczka:1994FoP,Pawlowski/Raczka:1995} the two physicists proposed   a ``Higgs free model for fundamental interactions'', as they  described it.  This proposal is formulated in a strictly conformal setting. Although it is very interesting in itself, we cannot discuss it here in more detail.

\subsection*{Mass generation by coupling to gravity: Drechsler and Tann}
 A view  close to Cheng's, establishing a  connection between gravity and electroweak  fields by  Weyl geometry,   was developed  a decade later   by Wolfgang Drechsler and Hanno Tann at Munich. Drechsler had been active for more than twenty years in differential geometric aspects of modern field theory.\footnote{For example \cite{Drechsler:Fibre_bundles}.}  
Tann  joined the activity  during his work on his PhD thesis \cite{Tann:Diss},  coming from a background interest in geometric properties of the de Broglie-Bohm interpretation of quantum mechanics. In their  joint work \cite{Drechsler/Tann}, as well as in their separate publications \cite{Tann:Diss,Drechsler:Higgs}  Weyl geometric structures are used in a coherent way, clearer than in  most of the other  physical papers cited up to now. 

They arrived, each one on his own,  at the full expression for  the (metrical) energy momentum tensor  of the scalar field, including terms which resulted from varying the scale invariant Hilbert-Einstein term (containing the factor $ \xi^{-1} $).\footnote{The terms with factor  $ \xi^{-1} $  had been introduced in an ``improved energy-momentum tensor'' by  \citeasnoun{Callan/Coleman/Jackiw} in a more ad-hoc way; cf.   \cite[(372)]{Tann:Diss}. }
 \beqarr   T_{\phi} &=&    D_{(\mu } \phi ^{\ast} D_{\nu )}\phi - \xi^{-1} D_{(\mu } D_{\nu )} |\phi |^2  
  \label{energy-momentum phi}  \\
& & \quad \quad  - g_{\mu \nu }
 \left( \frac{1}{2}D^{\lambda }\phi ^{\ast}D_{\lambda} \phi -   \xi^{-1} \, D^{\lambda} D_{\lambda}(\phi ^{\ast} \phi)  + V(\phi) \right)
\; . \nonumber
\eeqarr

 In their common paper, Drechsler and Tann introduced fermionic Dirac fields into the analysis of Weyl geometry  \cite{Drechsler/Tann}. Their  gravitational Lagrangian had the form
\beq  \mathcal{L}_{\text{grav}} =  \mathcal{L}_{R} +  \mathcal{L}_{R^2} \label{Drechsler/Tann Lagrangian}
\eeq
with $\mathcal{L}_{R}$  identical to Hung Cheng's (\ref{Hilbert term Weyl geometric}), in addition to $  \mathcal{L}_{\phi} $   (\ref{phi term Weyl geometric})  (with coefficients  $\beta = \frac{1}{6}$, $\varepsilon = -1$). A quadratic term, $ \mathcal{L}_{R^2} = \tilde{\alpha } R^2  \sqrt{|det\, g|}  $, in the (Weyl geometric) scalar curvature, was
 added.\footnote{In the appendix Drechsler and Tann showed that the squared Weyl geometric conformal curvature $C^2 = C_{\lambda \mu \nu \rho }C^{\lambda \mu \nu \rho }$ arises from the conformal cuvature of the Riemannian component $_g\hspace{-1pt}C^2$  by adding a scale curvature term: $C^2 = \, _g\hspace{-1pt}C^2 +\frac{3}{2} f_{\mu \nu }f^{\mu \nu }$ \cite[(A 54)]{Drechsler/Tann}. So one may wonder, why they did not replace the square term $ \mathcal{L}_{R^2}$ by the Weyl geometric conformal curvature term   $\mathcal{L}_{\text{conf}} = \tilde{\alpha } C^2  \sqrt{|det\, g|}$.  }

 For the development of a Weyl geometric theory of the Dirac  field, 
 Drechsler and Tann  introduced an adapted Lagrangian
\beq  \mathcal{L}_{\psi} =  \frac{i}{2} \left (\psi^{\ast}  \gamma ^{\mu}D_{\mu}\psi - D^{\ast}_{\mu}\psi^{\ast}  \gamma ^{\mu} \psi \right) + \gamma  |\Phi|  \psi^{\ast}\psi \, \label{Lagrangian Dirac field} \eeq  
with (scale invariant) coupling constant $ \gamma $ and Dirac matrices $\gamma^{\mu}$ with symmetric product $\frac{1}{2} \{ \gamma^{\mu},\gamma^{\nu} \} = g^{\mu \nu} \mathbf{1} $ \cite[(3.8)]{Drechsler/Tann}. Here the covariant derivative had to be  lifted to the spinor bundle, It  included already an additional  $U(1)$ electromagnetic potential $A=(A_{\mu})$
\beq  D_{\mu} \psi = \left( \partial _{\mu} + i \tilde{\Gamma} _{\mu} + w(\psi) \varphi_{\mu}  + \frac{i q}{\hbar c} A_{\mu}   \right) \psi \; ,  \eeq 
$q$ electric charge of the fermion field, $w(\psi) =-\frac{3}{2}$, $ \tilde{\Gamma} $ spin connection lifted from the Weylian affine connection.  This amounted to a (local) construction of a spin $\frac{1}{2}$ bunde. Assuming the underlying spacetime $M$ to be  spin, they worked in a Dirac spin bundle $\mathcal{D}$ over the Weylian manifold $(M, [(g,\varphi )]$. Its structure group was  $G= Spin(3,1) \times R^+ \times U(1) \cong Spin(3,1) \times \C^{\ast}$, where $ \C^{\ast}=\C \setminus {0}$.\footnote{One could then just as well consider a complex valued connection $z= (z_{\mu})$ with values $z_{\mu}= \varphi_{\mu} + \frac{i}{\hbar c} A_{\mu}$  in   $\C= Lie(\C^{\ast})$
and weight $W(\psi)= (-\frac{3}{2}, q)$. Then $D_{\mu} \psi= (\partial _{\mu} + \tilde{\Gamma}_{\mu} + W(\psi)z_{\mu})\psi$, presupposing an obvious convention for applying $ W(\psi) z$.  }

Drechsler and Tann considered (\ref{Lagrangian Dirac field}) 
as Lagrangian of a ``massless'' theory, because the masslike factor of the spinor field $\gamma  |\Phi|$ was not scale invariant.\footnote{This argument is possible, but not compelling.  $\gamma  |\Phi|$ has the correct scaling weight of mass and may be considered as such.} 
So they proposed to  proceed to a theory with masses by   introducing a ``scale symmetry breaking'' Lagrange term 
\beq  \mathcal{L}_{B} \sim \frac{R}{6} + (\frac{mc^2}{\hbar})^2 |\Phi|^2  \label{Lagrangian Weyl symmetry break}  \eeq
with fixed (non-scaling) $m$ \cite[sec. 4]{Drechsler/Tann}.\footnote{ So already in Tann's PhD dissertation.}

They did not associate such a transition from a seemingly ``massless'' theory to one with masses  to a hypothetical  ``phase transition''.
At the end of the paper they  commented:
 \begin{quote}
It is clear from the role the modulus of the scalar field plays in this theory (\ldots) that the scalar field with nonlinear selfcoupling is not a true matter field describing scalar particles. It is a universal field necessary to establish  a scale of length in a theory and should probably not be interpreted as a field having a particle interpretation. \cite[1050]{Drechsler/Tann}
\end{quote}
Their  interpretation of the scalar field  $\Phi$ was rather geometric than that of an ordinary quantum field; but their term
 (\ref{Lagrangian Weyl symmetry break}) looked ad-hoc to the uninitiated.\footnote{Note that one could just as well do without  (\ref{Lagrangian Weyl symmetry break}) and proceed with fully scale covariant masses -- compare last footnote. }  

\subsection*{Drechsler on mass acquirement of electroweak bosons}
Shortly after the common article appeared,  the senior author extended the investigation to  gravitationally coupled electroweak theory \cite{Drechsler:Higgs}.
Covariant derivatives were   lifted  as $\tilde{D}$ to the electroweak bundle. It included the additional connection components and coupling coefficients $ {g}$ and $ {g}'$ with regard to $SU(2)$ and $U(1)_Y$ like in Hung Cheng's work (\ref{ew covariant derivative}).
The Weyl geometric Lagrangian  could be generalized and transferred to the electroweak bundle  \cite[(2.29)]{Drechsler:Higgs},
\beq \mathcal{L}  =   \mathcal{L}_{\text{grav}}  + \mathcal{L}_{\Phi } + \mathcal{L}_{{\psi}} + \mathcal{L}_{\text{YM}} \, ,  \label{total Lagrangian ew-grav} \eeq
with contributions
like in (\ref{Drechsler/Tann Lagrangian}),  (\ref{phi term Weyl geometric}), (\ref{Lagrangian Dirac field}),  and (\ref{Cheng's ew interaction Lagrangian}) ({\em ew} terms only).
 Lagrangians for the fermion fields had to be rewritten similar to electromagnetic Dirac fields (\ref{Lagrangian Dirac field}) and were decomposed into the chiral left and right contributions. 

In principle, Drechsler's proposal coincided with Cheng's; but he proceeded with more care and with  more detailed explicit  constructions. 
He derived the  equations of motion with respect to all dynamical variables \cite[equs. (2.35) -- (2.41)]{Drechsler:Higgs} and calculated the energy-momentum tensors of all  fields ocurring in the Lagrangian. 

The symmetry reduction from the electroweak group $G_{\text{ew}}$ to the electromagnetic $U(1)_{\text{em}}$ could then be expressed    similar to the procedure in the  standard model. $SU(2)$ gauge freedom  allows  to chose a (local) trivialization of the electroweak bundle such that the $\Phi$  assumes the  form considered in the ordinary Higgs mechanism 
\beq \hat{ \Phi} \doteq  \left(  \begin{array} { c}   0  \\  \phi_o   \end{array} \right) \; ,
\eeq 
where $ \phi_o$ denotes a real valued field, and ``$\doteq$'' equality in a specific gauge. $ \hat{\Phi} $ has the isotropy group $U(1)$ considered as  $U(1)_{\text{em}}$. Therefore Drechsler called $ \hat{\Phi} $   the {\em electromagnetic gauge} of  $\Phi$.\footnote{In other parts of the literature  (e.g., the work of R\c{a}czka and Paw{\l}owski) it is  called ``unitary gauge'', cf. also  \cite{Flato/Raczka}. }

In two respects Drechsler  went beyond what had been done before: 
\begin{itemize}
\item He {\em  reconsidered the standard interpretation} of  symmetry breaking by the Higgs mechanism \cite[1345f.]{Drechsler:Higgs}. 
\item  And he  {\em calculated} the consequences of nonvanishing electroweak curvature components for  the 
{\em energy-momentum tensor } of the  scalar field    $\hat{\Phi}$ \cite[1353ff.]{Drechsler:Higgs}. 
\end{itemize}

With regard to the first point, he  made clear that he saw nothing compelling in the  interpretation of  symmetry reduction as ``spontaneous symmetry breaking due to a nonvanishing  vacuum expectation value of the scalar field''  \cite[1345]{Drechsler:Higgs}. He  analyzed the situation and came to the conclusion that 
the transition from our ${\Phi}$  to $\hat{\Phi}$  is to be regarded as a  ``choice of coordinates'' for the representation of the scalar field in the theory and has, in the first place, nothing to do with a ``vacuum expectation value'' of this field.\footnote{Mathematically spoken,  it is a change of trivialization of the $SU(2)\times U(1)$-bundle.}
\begin{quote}
 \ldots This choice is actually not a breaking of the orginal $\tilde{G}$ gauge symmetry [our $G_{\text{ew}}$, E.S.] but a different realization of it. (ibid.)
\end{quote}
He compared the stabilizer $U(1)_{\text{em}}$ of $\hat{\Phi}$ with the ``Wigner rotations'' in the study of the representations of the Poincar\'e group.
With regard to the second point, the energy-momentum tensor of the scalar field could be calculated roughly like 
 in the simpler case of a complex scalar field, (\ref{energy-momentum phi}).   Different to what one knew from the pseudo-Riemannian case,  the covariant derivatives $D_{\mu}\Phi$ etc. in (\ref{energy-momentum phi})  were then  dependent on  scale or $U(1)_{\text{em}}$ curvature. 

After breaking the Weyl symmetry by a Lagrangian of form (\ref{Lagrangian Weyl symmetry break})  (ibid. sec. 3), Drechsler 
 calculated the curvature contributions induced by the Yang-Mills potentials of the $ew$ group and its consequences for the energy-momentum tensor $ T_{\phi}$ of the scalar  field.
 Typical contributions to  components of $ T_{\phi}$ had the form of mass terms 
\beq m_W^2 W^{+  \, \ast}_{\mu} W^{- \,  \mu}, \; m_Z  Z_{\mu}^{\ast} Z^{\mu}  \, , \quad  \mbox{with} \quad  m_W^2 = \frac{1}{4} g^2 |\phi_o|^2  , \;
m_Z^2 = \frac{1}{4} g_o^2 |\phi_o|^2  \, , \label{Drechsler's mass terms}
\eeq
 $g_o^2 = g^2 +  g'^2 $, for the bosonic fields $W^{\pm}, Z$ corresponding to the generators $\tau _{\pm}, \tau _o$ of the electroweak group,  \cite[1353ff.]{Drechsler:Higgs}.\footnote{$W^{\pm}_{\mu} = \frac{1}{\sqrt{2}} (W^1_{\mu} \mp i W^2 _{\mu})$, $Z_{\mu}= \cos \Theta \,  W^3_{\mu} - \sin \Theta \,  B_ {\mu} $. } 
 They are identical with the mass expressions for the $W$ and $Z$ bosons in  conventional electroweak theory.
According to  Drechsler, the  terms (\ref{Drechsler's mass terms}) in $ T_{\phi}$ indicate  that 
the`` boson and fermion mass terms appear in the total energy-momentum tensor'' through the energy tensor  of the scalar field after ``breaking the Weyl 
symmetry''.\footnote{One has to be careful, however.  Things become more complicated if one considers the trace. In fact,  
 $ tr \, T_{\phi}$ contains a  mass terms of the Dirac field of form $\gamma | \phi_o | \hat{\psi}^{\ast} \hat{\psi}$, with $\gamma$ coupling constant of the Yukawa term
 ($\hat{\psi  }\;$  indicating electromagnetic gauge). That should be interesting for workers in the field.   One of the obstacles for  making quantum  matter fields compatible with classical gravity is  the vanishing of  $tr \, T_{\psi}$, in contrast to the (nonvanishing) trace of the energy momentum tensor of  classical matter.  
 Drechsler's analysis may indicate   a way out of this impasse. Warning: The mass-like expressions for $W$ and $Z$  in (\ref{Drechsler's mass terms}) cancel in $tr\, T_{\phi}$ \cite[equ. (3.55)]{Drechsler:Higgs} like in the energy-momentum tensor of the $W$ and $Z$ fields themselves.
In this sense, the mass terms of fermions and those of  electroweak bosons behave  differently with regard to the energy momentum tensor $T_{\phi}$.}
Inasmuch as the scalar field can be considered as extension of the gravitational structure of spacetime, the scale covariant theory of mass acquirement indicates a way to {\em mass generation}  by coupling to the gravitational structure. In any case, one has to keep in mind that the scalar field ``\ldots should probably not be interpreted as a field having a particle interpretation'' \cite[1050]{Drechsler/Tann}. 

Such a type of mass generation would have remarkable observable consequences in the LHC regime. The decay channels involving the standard Higgs particle would be completely suppressed.\footnote{A calculation of radiative corrections in the closely related conformal approach is presented in \cite{Pawlowski/Raczka:1995_II}; comparison  with  \cite{Kniehl/Sirlin:2000} might  be informative for experts.} If the LHC experiments turn out as a giant null-experiment with regard to chasing the Higgs particle, the  scale covariant scalar field should run up as a serious theoretical alternative to the Higgs mechanism.

\section{Scale covariance in recent cosmology   \label{section cosmology}}
\subsection*{Recent uses of Weyl geometry in cosmology}
 Already  early in the 1990s   Rosen and Israelit  studied different possibilities for  ``generating''  dark matter in 
Dirac's modified Weyl geometric framework  \cite{Israelit/Rosen:Dark_matter_1992,%
Israelit/Rosen:Dark_matter_1993,Israelit/Rosen:Dark_matter_1995,Israelit/Rosen:Particle}. They presupposed a  non-integrable scale connection leading to a spin 1 boson field which satisfied a scale covariant Proca equation, like in  Utiyamah's, Smolin's and Hung Cheng's papers. The authors called the new hypothetical bosons {\em Weylons} and proposed a crucial role for them in the constitution of  dark matter. 
In  recent years M. Israelit has developed ideas, how  matter may even have been ``generated from geometry'' in the very early universe \cite{Israelit:Matter_Creation_2002}\footnote{Compare with H. Fahr's proposal in his contribution to this volume.} and added a ``quintessence'' model in the framework of the Weyl-Dirac geometry \cite{Israelit:Quintessence_2002}.   Not all of it is convincing; but here is not the place to go into details.
 
 Weyl geometry has been reconsidered also by other authors as a possibility to relax the structural restrictions of Einsteinian gravity in a natural and, in a sense, minimal way. That happened indepently at several places in the world, at Tehran, Beijing, Santa Clara, Wuppertal, Atlanta, and perhaps elsewhere.  Some of these attempts built upon the Rosen/Israelit tradition of Weyl-Dirac geometry, others linked  to the field theoretic usage of Weyl geometry in  the standard model of elementary particle physics, or  to the Weyl geometric interpretation of the de Broglie-Bohm quantum potential. 

Entering the new millenium, our selective report will definitely leave the historical terrain in the proper sense. By pure convention I  consider work after  the watershed of the year 2000   as ``present''. It may, or may not, become object of historical research in some more or less distant future.  
The remaining section  concentrates  on those contributions of scalar fields or scale covariant aspects in cosmology  which relate directly or indirectly  to more basic aspects of Weyl geometry.

The scalar-tensor approach to gravity in the sense of Brans-Dicke theory has been studied all over the world. Among those active in this field, Israel Quiros from the university at Santa Clara (Cuba)   realized that the ``transformations of units'' in the sense of Brans-Dicke finds its most consequent expression in Weyl geometry. In some papers around the turn of the millenium he argued in  this sense  \cite{Quiros:2000a,Quiros:2000b};  but his main work  remains in more mainstream field theory and cosmology.

A turn of longer endurance toward Weyl geometry  was taken by M. Golshani,  F. and A. Shojai, from  the Tehran theoretical physics community. Their interest  stabilized when they studied  the link  between Brans-Dicke type  scalar tensor theory and de Broglie-Bohm quantum mechanics. About 2003 (perhaps  during their stay at the MPI for gravitation research Golm/Potsdam) the Shojais realized that Weyl geometry can  be used as a  a unifying frame for such an enterprise \cite{Shojai/Shojai:2002}. It seems that their retake of Weyl geometry  may have influenced other colleagues of the local physics community, who started   to analyze astrophysical questions by Weyl geometric methods \cite{Moyasari_ea:2004,Mirabotalebi_ea:2008}.  The  Weyl geometric background   knowledge  of the Tehran group  was shaped by Weyl-Dirac theory and the Rosen/Israelit tradition, supplemented by the analysis of   \citeasnoun{EPS} and  of \citeasnoun{Wheeler:1990}. The latter had explored the relation between quantum physics and (Weyl) geometry already back in the 1990s.

E. Scholz, a historian of mathematics at the Mathematics Department  of  Wuppertal University  (Germany)  started studying   Weyl geometry in present cosmology shortly after attending a conference  on history of geometry at Paris in September 2001.\footnote{On this conference P. Cartier, a  protagonist of the second generation of the Bourbaki group who has been interested in mathematical physics all his life, gave an enthusiastic  talk on the importance of Weyl's scale connection for understanding cosmological redshift \cite{Cartier:Cosmology}. Scholz was struck by this talk, because he had tried  to win over  physicists for such an idea in  the early 1990s, of course without any success.}
 Coming from a background in mathematics and its  history,  it took some time before he  got  acquainted with the more recent Weyl geometric  tradition in theoretical physics.  After he ``detected''   the work of  Drechsler and Tann on Weyl geometric methods in field theory  in late 2004, it  became a clue for his entering the physics discourse  in field physics \cite{Scholz:model_building_2005,Scholz:FoP2009}. 

C. Castro had become  acquainted with Weyl geometry in physics already in the early 1990s while being at Austin/Texas. At that time a proposal by E. Santamato's  to use Weyl geometry for a geometrization of de Broglie-Bohm quantum mechanics stood at the center of his interest \cite{Castro:1992}.\footnote{\cite{Santamato:KG,Santamato:WeylSpace} } After the turn of the millenium, then working  at Atlanta (Centre for Theoretical Studies of Physical Systems), he took up the Weylian thread again, now with the guiding questions, how Weyl's scale geometry may be  used for gaining  a deeper understanding of dark energy and, perhaps, the Pioneer anomaly \cite{Castro:FoP2007,Castro:2009}.

A completely different road towards  Weyl geometry was opened for Chinese theoretical physicists Hao Wei, 
Rong-Gen Cai, and others   by a talk of Hung Cheng, given in July 2004  at the Institute for Theoretical Physics of the Chinese Academy of Science, 
Beijing.\footnote{\cite[Acknowledgments]{Wei_ea:2007}} 
It was natural for them to take the ```Cheng-Weyl vector field'' (i.e., the Weylian scale connection with massive boson studied by Cheng in the late 1980s) and Cheng's view of the standard model of elementary particle physics as their starting point \cite{Wu:2004,Wei_ea:2007}.

So far only  groups or persons have been mentioned, who  contributed explicitly   to the present revitalization of  Weylian scale geometry. Other protagonists whose work plays a role for this question will enter this  section, even if they do not care about  links  to  Weyl geometry.

\subsection*{Mannheim's conformal cosmology}
A  striking analysis of  certain aspects of recent cosmology (the {\em dark} ones, dark matter and dark energy) was given by Philip Mannheim and Demosthenes Kazanas. In the 1980s the two physicists analyzed the  ``flat'' shape of galaxy rotation curves (graphs of the rotation velocity $v$ of stars in dependence of the distance $r$ from the center of the galaxy). From a certain distance close to a characteristic length of the galaxy ($2.2\, r_o$ with $r_o$ the ``optical disc length scale'') $v$ is greater than expected by Newtonian mechanics, such that the spiral should have flewn apart unless unseen (``dark'') matter enhancing  gravitational binding or a modification of Newtonian/Schwarzschild gravity were assumed. While the  majority of astrophysicists and astronomers assumed the first hypothesis, Mannheim looked for possible explanations along the second line.   In \cite{Mannheim/Kazanas:1989} a theoretical explanation of the flat rotation curves was given, based on a conformal approach to gravity. During the following years the approach was deepened and extended to the question of ``dark energy''.

In fact Mannheim and Kazanas found that, in  the conformal theory, a static spherically symmetric matter distribution could be described by the solution of a fourth order Poisson equation
\beq  \nabla ^4 B(r) = f(r) \eeq
  with a typical coefficient  $B(r)$ 
\[ B \sim  - g_{oo}=  g_{rr}^{-1} \]
of a metric $ds^2  = g_{oo}dt^2 -  g_{rr}dr^2 - r^2 d\Omega ^2$ (up to a conformal factor). 
  The r.h.s. of the Poisson equation, $f(r)$, depended on the mass distribution, e.g., in a spiral galaxy.

A general solution turned out to be of the form
\beq g_{oo}=    1 - \beta (2-3\beta \gamma ) r^{-1} - 3 \beta \gamma + \gamma r - \kappa r^2 \,  \eeq
with constants $\beta ,\gamma ,\kappa $. Here $\beta $ depends  on the mass and its distribution in the galaxy. For galaxies   Mannheim arrived at such small values for  $\beta $ and $\gamma $ that the $\beta \gamma $ terms could be numerically neglected, as could the $r^2$ term. In this case the  $\beta $ term  took on the form of a Schwarzschild solution of the Einstein equation with Schwarzschild radius $r_S = \beta $. The classical potential  was, however,  modified by a  term linear in $r$, in addition to the classical Newton potential,
\beq   V(r) = - \frac{\beta }{r} + \frac{\gamma }{2} r   \,  \eeq 
and a corresponding velocity of generalized Kepler orbits
\beq v(r) \sim V(r) \, .  \label{Mannheim's rotation curves}
\eeq 
The dynamics of such a potential agreed well with the  data of galactic rotation curves \cite{Mannheim/Kazanas:1989,Mannheim:1994}. 

The result for 11 galaxies with different behaviour of rotation curves led to  surprisingly good fit. $\gamma$ was basically independent of the galaxy. It  had  a cosmological order of magnitude, $\gamma \approx 3 \cdot 10^{-30}\, cm^{-1} \approx 0.04 H_1$. 
 Mannheim considered this ``an intriguing fact which suggests a possible cosmological origin for $\gamma $'' \cite[498]{Mannheim:1994}. In his view, it represented a kind of (weak) Machian type influence of very distant masses on the potential, which could be neglected close to stars and at the center of the bulge  of galaxies. It came to bear only at the galactic  periphery and beyond.

With respect to the dark energy problematics, Mannheim chose a peculiar perspective. In the special case of conformally flat models, like Robertson-Walker geometries, he decided to consider the scale invariant Hilbert-Einstein Lagrangian $-\frac{1}{12}|\phi |^2 R \sqrt{|det\, g|}  $ as part of the {\em matter} Lagrangian. Due to his sign choices, he arrived at a version of the Einstein  equation with {\em inverted } sign and interpreted it as a  kind of ``repulsive gravity'' which he claims to operate on cosmic scales, in addition to ``attractive gravity'' on smaller scales  indicated by the conformally modified Schwarzschild solution. In his eyes, such a repulsive gravity might step into the place of the dark energy of the cosmological constant term of standard gravity \cite[729]{Mannheim:2000}.

In spite of such a grave difference to Einstein gravity, Mannheim does not consider his conformal approach to   disagree  with the standard model of cosmology and its accelerated expansion. He rather believes that his approach may lead to a more satisfying explanation of the expansion dynamics. In his view, ``repulsive gravity'' would take over the role of dark energy. Moreover he expects that it may shed new light on the  initial singularity and, perhaps, also on the black hole singularities inside galaxies.

\subsection*{Scalar fields and dark energy: Kim and Castro}
Other authors started to analyze dark energy by a scalar field approach.
 A  remarkable contribution to this topic comes from  Hongsu Kim (Seoul). He uses 
a classical Jordan-Brans-Dicke field, $\chi =  |\phi|^2$ with JBD parameter $\omega $ like in (\ref{Lagrangian JBD}) and shows  that, under certain assumptions, it may lead to a phase of {\em linear expansion} in   ``late time'' development of the cosmos, i.e. long after inflation but long before ``today'' \cite{Kim:dark_energy} . He  proposes to consider a transitional phase  between a matter dominated phase with decelerated expansion (decelerated because of dominance of gravitational attraction over the repulsive vacuum energy) and an  accelerated expansion dominated by  vacuum energy. For the JBD field he uses the funny terminology of 
{\em k-essence},  in distinction to {\em quintessence} which has been introduced for scalar fields without coupling to gravity.

Assuming a spatially flat Robertson-Walker spacetime with  warp (scale) function $a(t)$, without ordinary matter and without cosmological constant, Kim starts from the modified Friedmann  and scalar field equations:
\beqa
\left( \frac{\dot{a}}{a} \right)^2 &=& \frac{\omega }{6}\left(\frac{\dot{\chi}}{\chi } \right)^2 -  \frac{\dot{a}}{a} \left(\frac{\dot{\chi}}{\chi } \right) \\
\ddot{\chi } + 3  \frac{\dot{a}}{a} \chi &=& 0 \, 
\eeqa
($t$  time parameter,    $\chi$  scalar field,  $\omega $  JBD parameter). 

By an  {\em Ansatz} of the form $a(t) \sim (at +b)^{\alpha }, \; \chi (t) \sim a(t)^n$ he evaluates the energy stress tensor of $\chi$. The  divergence condition $\nabla ^{\nu}T_{\mu \nu }^{(\chi )}=0$ implies the restriction $\omega =-\frac{3}{2}\, \; n = -2$. This leads to a solution with linear warp function
$  a(t) \sim t, $ 
and 
$\rho _{\chi }\sim a(t)^{-2}\, , \qquad p_{\chi } = - \frac{1}{3}  \rho _{\chi } $ 
for the  energy density $\rho _{\chi }$ and pressure $p_{\chi }$ of the scalar field.  According to the author, this solution indicates a kind of dynamically neutral intermediate  state of the underlying universe model between deceleration and acceleration, which does not arise in the standard model \cite[equs. (8), (13), (19)]{Kim:dark_energy}.

Kim goes on to consider a Lagrangean with JBD field term (\ref{Lagrangian JBD}), including ordinary matter $m$ and cosmological constant $\Lambda $.
He assumes that the ``late'' time evolution of the cosmos consists of three phases: 
\begin{itemize}
\item at the beginning of ``late time'' matter  dominates the other two terms (deceleration), 
\item in an  intermediate stage  matter density has been diluted sufficiently far so that the scalar field dominates the dynamics (linear expansion,  $m$ and $\Lambda $ negligible), 
\item after further dilution of $\rho _{\chi }$,  the vacuum term $\Lambda $ takes over and   dominates the  evolution of the model   (acceleration). 
\end{itemize}

In the first phase $\chi$  ``mixes'' with (ordinary) matter. Kim considers this as a state of dark matter. In the third phase  $\chi$ seems to mix with the   vacuum term. The scalar field is then  assimilated to  dark energy. In this  sense, so Kim argues,
 JBD theory offers ``a unified model'' for dark matter and dark energy.

 Kim's analysis of the contribution of JBD fields to cosmological dynamics is highly interesting; but it has two empirical drawbacks. Firstly and noticed by himself, the present standard model does not know of any linear phase of expansion.\footnote{This is no great handicap because the  paradigm of standard cosmology does not allow such a phase and therefore could not  ``see'' it, even if it existed. Only extreme observational results could {\em enforce} such a phase  onto the standard view and would  then break it up.}
Much worse, although not discussed by Kim, is the empirical restriction for the JBD parameter $\omega > 10^3 $, following from the comparison of post-Newtonian parametrized gravity theories with the data of  high precision observations.\footnote{Cf. the contribution of C. Will, this volume.}
 This observational result stands in  grave contradiction to Kim's theoretical value 
$\omega =-\frac{3}{2}$.\footnote{One should keep in mind, however,  that the observational constraint for $\omega $  refers to {\em Jordan frame},  which  Kim  presupposes  at the moment. In the Weyl geometric reading  of  JBD theory the  Einstein frame would appear more appropriate  for observational purposes (see below, under equ. (\ref{observables})). }

Some of our authors assume 
\beq D_{\mu}\phi = 0  \; , \label{trivial scale connection}
\eeq 
for the sake of simplicity  \cite{Kao:1990}, others by their background in  de Broglie-Bohm theory \cite[equ. (14)]{Santamato:WeylSpace}. 
In the context of scalar tensor theory, however, this condition 
  introduces a dynamical overdetermination and makes Weyl geometry essentially redundant. This  can  be seen by comparison with the JBD approach. In Riemann gauge,  condition (\ref{trivial scale connection}) implies
$\partial_{\mu} \phi = 0$  and thus $|\phi|= const$.  In terms of  of Brans-Dicke theory one arrives at  a constant scalar field in Jordan frame (something like a {\em contradictio in adjecto})!  Then the latter is identical with  the Einstein frame, and the whole JBD theory  becomes trivial. 
In the sequel I shall call assumption of (\ref{trivial scale connection}) a {\em trivial Weylianization} of Einstein gravity.

Castro arrives at condition (\ref{trivial scale connection}) by 
varying the Weyl geometric Lagrangian  (\ref{Hilbert term Weyl geometric}), (\ref{phi term Weyl geometric}) 
\[  \mathcal{L} =   \mathcal{L}_{R}  + \mathcal{L}_{\Phi } + \mathcal{L}_m \]
not only with respect to $g$,  $\phi$ and  the matter variables, but also with respect to $\varphi$ \cite[equ. (3.21) ff.]{Castro:FoP2007}, \cite[equ. (10)]{Castro:2009}. 
As we have just seen,  $\varphi$ is  no independent dynamical variable in the scalar field approach to {\em integrable} Weyl geometry. It is nothing but the ``other side'' of the scalar field which has  a  dynamical (Klein-Gordon) equation of its own.\footnote{In order to avoid misunderstanding, it has to be added that $\varphi$ is, of course, a dynamical variable in the nonintegrable case. Then one has to introduce a scale curvature term into the Lagrangian (usually of Yang-Mills type, $\frac{1}{4} f_{\mu \nu }f^{\mu \nu }$ like in the work of Dirac, Smolin (\ref{Smolin's Lagrangian}), Hung Cheng and Drechsler/Tann. 
The integrable case arises  from constraining conditions which can be expressed in the action by a system of (antisymmetric) Lagrange multipliers $\lambda ^{\mu \nu }$   for  scale curvature $f=(f_{\mu \nu })$, 
 \beq \mathcal{L}_{f}= \left( \frac{1}{4} f_{\mu \nu }f^{\mu \nu } + \frac{1}{2} \lambda ^{\mu \nu }   f_{\mu \nu } \right) \sqrt{|det \, g|}\; .
\eeq 
The variation considered by Castro then acquires  an additional term from the derivatives of the Lagrange multipliers, and the result  looks something  like $D_{\mu}\phi = D_{\nu} \lambda ^{\nu}_{\mu}$.
Variation with respect to the multipliers $\lambda ^{\mu \nu }$  gives the integrability constraint
$ f_{\mu \nu }= 0 $. }
With some additional assumptions on a de Sitter solution he ``derives''  a constant vacuum energy of the ``right order of magnitude''.

\subsection*{Taking rescaling  seriously: Scholz, Masreliez}
In JBD theory the physical consequences of conformal rescaling have long been discussed, with a tendency toward shifting from Jordan frame to Einstein frame as the preferred one for physical observations \cite{Faraoni_ea:Transformations}. A similar discussion in Weyl geometry is nearly lacking. One of the exceptions is a passage in \cite{Scholz:FoP2009}. This paper  also presents  a 
 quite different direction for a Weyl geometric analysis of  questions relating to  dark energy.   Scholz,  starts from  Drechsler's and Tann's work but simplifies the Lagrangian for a  classical approch to cosmology.

 He abstracts from the particle fields of the SMEP and plugs in a classical matter term of an ideal fluid defined with respect to a timelike unit vector field $X=(X^{\mu})$ like in \cite[69f]{Hawking/Ellis}. These have to be  adapted to the Weyl geometric approach by ascribing  proper  scale weights, $w(\rho_m) =-4$ and $w(X)= - 1$,  to  matter energy density $\rho_m$ and to the unit field $X$.\footnote{Such an assumption is completely natural from   ``transformation of units'' view of scale transformations. \label{footnote transformation of units}}
 From  the results of Ehlers/Pirani/Schild and 
Audretsch/G\"ahler/Straumann Scholz draws the consequence of postulating integrability of the Weylian scale connection, $d \varphi = 0$.\footnote{\cite{EPS,AGS}}
Abstracting from radiation, in a first approach, he  arrives at a Weyl geometric (scale invariant) version of the Einstein equation,
\beq Ric - \frac{1}{2}R g =  (\xi |\phi|)^{-2}(T^{(m)} + T^{(\phi)})
\eeq
with  matter term $T_{\mu \nu }^{(m)}= (\rho_m + pm)X_{\mu }X_{\nu } + p_m g_{\mu \nu }$ well known from ideal fluids, but here appearing with scale covariant matter energy density and pressure $p_m$. The  energy stress tensor of the scalar field $T^{(\phi)}$ has been computed by Drechsler/Tann (\ref{energy-momentum phi}). $\xi$ is a coefficient regulating the relative coupling strengths of $\phi$ to the scale invariant Hilbert-Einstein term in comparison to its kinetic term, like in (\ref{Lagrangian JBD Fujii/Maeda}). In this paper Scholz  follows the prescription from conformal gravity that  $\xi$ should be set to $ \frac{1}{6} $ for spacetime dimension $n=4$, although in the Weyl geometric frame scale invariance of the Lagrangian holds  for any value of $\xi$; and $\xi$ is an adaptable parameter, at least at the outset.\footnote{The Weyl geometric Einstein-Hilbert Lagrangian $\xi |\phi|^2R\sqrt{|det\, g|}$ and the kinetic term $D^{\mu}\phi^{\ast}D_{\mu}\phi \sqrt{|det\, g|}$ are  indepently scale invariant. In the conformal theory they are scale invariant only as a joint package after the  choice of $\xi = \frac{n-2}{4(n-1)} $.} 

Close to Weyl's ``calibration by adaptation''  Scholz introduces ``scale invariant magnitudes''  $\hat{Y}$  of any scale covariant quantity $Y$ (scalar, vectorial, tensorial etc.)  at a point $p$ by\footnote{Cf.  footnote \ref{footnote transformation of units}.}
\beq \hat{Y}(p) =  |\phi(p)|^{w}\, Y(p) \, ,  \qquad w:= w(Y) \, . \label{observables}
\eeq 
Observable quantities can be evaluated in any scale gauge; $\hat{Y}$ is gauge invariantly defined. In this sense all  gauges have equal status.
 But the  determination is easiest in scalar field gauge with $|\phi| = const$, corresponding to ``Einstein frame'' in the terminology of JBD theory, because in this gauge $\hat{Y}(p)\doteq Y(p)$.\footnote{Cf.  Utiyama's and Israelit's view of $\phi$ as a ``measure field''.}

 Robertson-Walker cosmologies  arise here from the usual assumption of a global  foliation with homogenous and isotropic spacelike folia 
orthogonal to the timelike vector field $X$, which is now  identified with an  observer field specified by the flow.

Scholz hints at the striking property that cosmological redshift $z$ of a signal emitted at a point $p_1$ and observed at $p_2$ (with respect to the observer field) after following a null geodesic $\gamma$ is scale invariant,\footnote{The invariance of $z$  is due to the natural scaling of the timelike field, the invariance of  (null and other)  geodesics and the scaling of the metric.} 
\beq z+1 = \frac{g_{p_1}(\gamma '(p_1),X(p_1))}{g_{p_2}(\gamma '(p_2),X(p_2))}  \, .
\eeq 
Thus the Weyl geometric look at Brans-Dicke fields shows  two things: Firstly,  the rescaling to scalar field gauge (Einstein frame) is ``natural'' in the sense of giving direct access to observable magnitudes, although not the only one in which observable magnitudes can be calculated. Second, cosmological redshift is not exclusively bound to ``expansion'' (warping by $a(t)$); it may just as well depend on the scale connection or  the scalar field. In $|\phi|$-gauge it is an effect composed by  a warp function contribution  and by   Weyl's scale transfer function 
\beq  \lambda (p_1,p_2) = e^{\int _1^{2} \varphi(\gamma ')} \,  \label{scale transfer} \eeq
 ($\gamma $ path between points 1 and 2).

In certain cases the residual expansion vanishes and cosmological redshift is exclusively given by the scale transfer, $z+1 = \lambda$. In this case  $z$ can just as well be attributed to the scalar field (in Riemann gauge/Jordan frame) as to  the Weylian scale connection linked to it (in scalar field gauge/Einstein frame). {\em Expansion no longer appears  as a physically real effect of cosmology.}
The warp function $a(t)$ is  reduced to an auxiliary role in the mathematics of Robertson-Walker-Weyl geometries.

The final aim of Scholz' investigation is a study of those Robertson-Walker-Weyl geometries in which such an extreme reduction of cosmological redshift to the scale connection (respectively the scalar field) happens. In honour of the inventor of the mathematical framework used by him  he  calls them ``Weyl universes''. Summarily stated a {\em Weyl universe} is given, in scalar field gauge, by a static spacetime geometry with Riemannian component $g$ of the Weylian metric with spatial slices of constant sectional curvature $\kappa = k a^{-2}$, $k= 0, \pm 1$  ($a$ ``radius of curvature''),
\beq   d{s}^2 = - d \tau^2 +   \frac{dr^2}{1-\kappa\, r^2} + r^2 ( d\Theta ^2 +  r^2 \sin^2\Theta \, d\phi^2 )  \, .
\eeq
The scale connection is time homogeneous with  only nonvanishing component  $H$ (constant) in time direction,
\beq \varphi = (H, 0,0,0) \, .
\eeq 
The scalar field is,  of course, also constant,
\beq |\phi|^2  = \xi^{-1} \frac{c^4}{8 \pi G} = 6 \frac{c^4}{8 \pi G} \, .
\eeq 
 For the transition to Riemann gauge $\tilde{g} = \Omega ^2g, \;  \tilde{\varphi}= 0$ (Jordan frame) the Weylian scale transfer function is used,   $  \Omega (t) = e^{Ht}$. That leads to the metric,
 \beq d{\tilde{s}}^2 = e^{2Ht}\left( - d \tau^2 +   \frac{dr^2}{1-\kappa\, r^2} + r^2 ( d\Theta ^2 +  r^2 \sin^2\Theta \, d\phi^2 )  \right)\, , \label{SEC metric}
 \eeq 
 called {\em scale expanding cosmos} by J. Masreliez (see below), and to an exponentially decaying scalar field $\tilde{\phi}(t)= e^{-Ht}$. 

A change of time coordinate, $\tau =e ^{Ht} $ shows that the ``scale expanding cosmos'' (\ref{SEC metric}) is nothing but a {\em linearly expanding Robertson-Walker-Weyl model}. It has an inversely decaying scalar field which influences observable quantities in the sense of (\ref{observables}),
 \beqarr \qquad \quad  d{\tilde{s}}^2 &=& - d \tau^2 + (H\tau)^2 a^2 \left(  \frac{dr^2}{1- k\, r^2} + r^2 ( d\Theta ^2 +  r^2 \sin^2\Theta \, d\phi^2 )  \right) , \\
  \tilde{\phi}(\tau) &=& (H \tau)^{-1}  \,  , \qquad    \tilde{\varphi} = 0 \, .
 \eeqarr 
 The observable  Hubble parameter is   $\hat{H} = H$.\footnote{The Hubble coefficient $\tilde{H}$  is measured as a (reciprocal) energy reduction of photons over  distance (resp. running time) and has  dimension  inverse  length. It is  a magnitude of scale weight $w(\tilde{H})=-1$. The Hubble parameter of ordinary Robertson-Walker theory is $H(\tau) = \frac{\dot{a}}{a}  = \tau ^{-1}$. Its  Weyl geometric observable magnitude (\ref{observables}) is thus $\hat{H}(\tau) = |\phi (\tau)|^{-1} H(\tau)= H$. }
 
Close to the end  of \cite{Scholz:FoP2009} the author  investigates conditions under which  the energy stress tensor of the  scalar field stabilizes  an Einstein-Weyl universe (i.e. one with positive sectional curvature $\kappa >0$). He comes to the conclusion that this may happen if a a certain relation between $H$, mass energy density $\rho_m$  and the coefficient $\lambda _4$  of the fourth order term of the scalar field potential ($V(\phi) = \lambda _4 |\phi|^4$) is satisfied \cite[64]{Scholz:FoP2009}. This  condition seems neither particularly natural nor theoretically impossible. Although one may reasonably doubt that this is the end of the story, it demonstrates 
 the existence of  unexpected possibilites for the energy stress tensor of the scalar field in the Weyl geometric approach. In fact, they are excluded  in  Einstein gravity by the singularity theorems of Penrose and 
Hawkins.\footnote{Kim's linearly expanding phase of a Robertson-Walker model with JBD field sheds light on Scholz' approach, and vice versa. Embedded in a Weyl geometric context, Kim's model is nothing but a Minkowski-Weyl universe with flat spatial slices 
($\kappa =0 $). If Kim's analysis of the pure scalar field dynamics without potential is correct, it tranfers to the Weyl geometric context after adaptation of the parameter.}

Conceptually  Weyl universes are closely connected to the theory of  {\em scale expanding cosmos} (SEC),  proposed by J. Masreliez \cite{Masreliez:Apeiron2004,Masreliez:Apeiron2004b}. Masreliez works with a metric  like in (\ref{SEC metric}) and tries to rebuild more or less the whole of cosmology.
He doubts the reality of cosmic expansion from a physicists point of view\footnote{At the moment this is  a minority position among physicists, best expressed by correspondents  of the {\em Alternative Cosmology Group},  www.cosmology.info. } and  argues for ``physical equivalence'' of the scale expanding metric very much in the sense of scale co/invariant theory:
\begin{quote}
[S]cale expansion for flat or curved spacetimes does not alter physical relationships; scaled spacetimes are equivalent and scale invariance is a fundamental, universal, gauge invariance. \cite[104]{Masreliez:Apeiron2004}
\end{quote}
Masreliez calls upon scale invariant theory and, in the end, of Weyl geometry. His SEC model is basically nothing but a Weyl universe, considered in Riemann gauge. He generally prefers a flat SEC, or in our terms, a Minkowski-Weyl universe. But he does not realize that Weyl geometry might be helpful for his enterprise.  We cannot enter into details here, because not all of Masreliez' explanations are   mature; many seem not  particularly clear to the reporting author.

\subsection*{Attempts for  understanding ``dark matter''}
Besides the  more widely noticed approach of modified Newtonian dynamics (MOND)  and Mannheim's conformal gravitation, some researchers  try to understand   the flat rotation curves of galaxies  by the contribution of  scalar fields to  the total energy  around galaxies and clusters.  Franz Schunck started with such research already in his Cologne PhD thesis under the direction of E. Mielke    \cite{Schunck:Diss,Schunck:1998}.  He continued this line of investigation in different constellations. Here we concentrate on joint work with   Mielke and  Burkhard Fuchs \cite{Mielke/Schunck_ea}.\footnote{Scale covariance of the scalar field does not  play a role in this research; but it should be just a question of time until this further step is undertaken.}
Already in the 1970s E. Mielke  investigated a Klein-Gordon field  $\phi$ with bicubic (order $6$) potential  
\beq V(\phi)= m^2 |\phi|^2(1- \beta |\phi|^4) \; , \qquad \beta |\phi| \leq 1 \, .  \label{Mielke's toy model}\eeq
He found  that  the corresponding nonlinear Klein-Gordon equation (with higher order
 self-interaction) allows for non-topolo\-gi\-cal soliton solutions \cite{Mielke:1978,Mielke:1979}. Our three authors could draw upon this result at least as a  ``toy model'', as they admit, for their investigation of scalar fields as a potential source for  {\em dark matter} \cite[44]{Mielke/Schunck_ea}. 
In an empirical investigation of their own Fuchs and Mielke  show that the observed rotation curves of 54 galaxies stand in good agreement with expectations derived even from their toy model \cite{Fuchs/Mielke:2004}.

The work of the three authors demonstrates that the scalar field approach to ``dark matter'' is worthwhile to pursue. It is concentrated and goes deeply  into  empirical evidence. But the gap between  scalar field halos and gravity remains still wide open. Some authors have  started recently  to tackle such questions  based on JBD or Weyl-Dirac theory.  

 Hongsu  Kim, whom we met already in the passage on dark energy, does  so in  \cite{Kim:2007} by the  Brans-Dicke approach. He uses a method developed already in the 1970s  to construct axis symmetric solutions of the scalar field equation and the  JBD version of Einstein equations with a singularity along the symmetry axis. He adds a  highly interesting discussion of the singularity of the modified metric along the axis ($\theta =0$). He remarks that this is a true singular direction, not only a coordinate singularity, and   estimates the velocity of timelike trajectories close to the axis. He finds them  to be close to the speed of light and arrives at the conclusion that the `` bizarre singularity at $\theta =0, \, \pi $ of the Schwarzschild-de Sitter-type solution in BD gravity theory can account for the relativistic bipolar outflows (twin jets) extending from the central region of active galactic nuclei (AGNs)'' \cite[24]{Kim:2007}.

If  this observation is right, even though only in principle, it will by far outweigh Kim's  rough estimation of rotation velocities. The acceleration of jet matter is an unsolved  riddle of astrophysics. It would be a great success for the approach, if a scalar field extension of Einstein gravity would be able to give a clue to this challenging phenomenon. 
Kim's   analysis is formulated in classical JBD theory. He does not even consider conformal transformations to the Einstein frame, but rather stays in the Jordan frame.  It will be interesting to see, what change it  will make to take up conformal rescaling and Weyl geometric methods for this question.

None of the scalar field or Weyl geometric attempts to explain ``dark matter'' can yet compete in precision with Mannheim's conformal gravity approach. But in the range of proposals (in particular of Fuchs/Miel\-ke/Schunck, Kim, and Castro) we have the seeds of an Ansatz for a promising research program. Whether it will lead to a solution of the problem remains, of course, to be seen.

\section{Discussion \label{section discussion}}
We have seen that JBD scalar tensor theory works in a Weyl geometric structure, although in most cases unknowingly.\footnote{The scale connection $\varphi$ in ``Einstein frame'' (scalar field gauge) is here usually hidden in partial derivative expressions equivalent to  $\varphi = - d \log \phi = - d \omega  $, where  $\phi(x) = e^{\omega (x)}$ is the scalar field in ``Jordan frame'' (Riemann gauge).} 
The analysis of Ehlers/Pirani/Schild  shows  that Weyl geometry is  deeply rooted in the basic structures of gravity.    A first step towards  founding  gravity upon quantum physical structures (flow lines of WKP approximation of Dirac or Klein-Gordon fields) rather than on classical particle paths has been made by 
Audretsch/G\"ahler/Strau\-mann.\footnote{\cite{EPS}, \cite{AGS}.} Deeper links (Feynman  path integral methods) or broader ones  (geometrization of  quantum potential)  to quantum physics have not been discussed in this 
article.\footnote{For the first question see \cite{Narlikar/Padmanabhan},  for the second \cite{Santamato:KG,Santamato:WeylSpace} and the succeeding literature. A recent survey on the last point is given in \cite{Carroll:grav_quantum_potential,Carroll:quantum_potential}.}
But already the  twin segments of theoretical physics considered  here, elementary particle physics and cosmology, show remarkable features of  Weyl geometry in recent and present physics.

There seem to be  intriguing perspectives for Weyl geometric scalar fields, at least on a theoretical level, for approaching the problem of  ``mass generation'' (as it is usually called) in particle physics.  Englert, Smolin,  Cheng, R\c{a}czka/Paw{\l}owski and  Drechsler/Tann have opened  a   view, not widely perceived among present physicists, of how a basic scalar field may participate in the mass generation of fermions and the electroweak bosons by coupling them to gravity. Their work has  been marginalized during  the rise of the standard model SMEP. But  time may be ripe for  reconsidering this  nearly forgotten Ansatz. 

From the Weyl geometric perspective,  Drechsler's and Tann's analysis  indicates most markedly   a possible  link between gravity and particle physics at an unexpectedly ``low'' energy level. This energy level will be  reached by  the  LHC experiments startin  November 2009. Already during the next  five to ten years we shall learn more about  whether the famous Higgs particle does indeed show up, at the end of the day, or  whether the scepticism of  scale covariant scalar field theoricians with regard to  a massive Higgs particle is  empirically supported in the long run.\footnote{Cf. \cite{Kniehl/Sirlin:2000}.}

We also have found interesting aspects of the analysis of the ``dark'' matter problem by scale covariant scalar fields  in the works of  Kim, Castro, Fuchs/Miel\-ke/Schunck and others.  Scalar field models of spherically or axially symmetric solutions of the slightly generalized Einstein equation (and the Klein-Gordon equation) show properties which open   promising perspectives for further investigations.  The authors of these researches come from different directions and start to dig a land which seems worth the trouble of farming it more deeply. From our point of view, it it should be  investigated whether   introducing  gauge transformations and  Weylian scale connections  into  Kim's approach, or considering not only ``trivial'' Weylianizations (as we have called it above) in Castro's approach,  helps to advance the understanding of the dark matter problem. 

Most remarkable seem  the structural possibilities opened up by Weyl geometry for  analyzing  how  the scalar field energy tensor introduces a repulsive element into (scalar field extended) gravity, usually considered as vacuum or ``dark''  energy.  The scalar field energy tensor allows for a much wider range of dynamical possibilities  than usually seen in the framework of classical Friedmann-Lemaitre models. Even a balanced (non-expanding) spacetime geometry appears to be dynamically possible and, under certain assumptions, even natural. It is interesting to see that Kim's linearly expanding Robertson-Walker model, Masreliez' scale expanding cosmos, and Scholz' Weyl universes characterize  one and the same class of spacetimes  in the framework of Weyl geometry. 

In addition,  Weyl geometric gravity theory, with modified scale invariant Hilbert-Einstein action coupled to the scalar field, sheds new light on cosmological redshift.  In this frame, the famous  expanding space explanation of the Hubble redshift appears  only as one possible perspective among others. From a theoretical point of view, it  even need  not be considered as  the most convincing one.\footnote{Remember Weyl's expectation that some day a ``more physical'' explanation will probably  take the place of the ``space kinematical''  description of cosmological redshift. He considered the latter as of    heuristic value only, due to its mathematical simplicity \cite{Weyl:Redshift}. }
With such questions we  enter a terrain which physicists usually consider as morass; but  Weyl geometry gives these investigations a safe  conceptual framework.

 Also here, like in the case of the Higgs mechanism, we have increasing observational evidence. It will contribute  either to dissolving the standard wisdom or harden it against  theoretically motivated scepticism. 
In the case of cosmological redshift, we look forward with great interest to more data on metallicity in galaxies and quasars. At present the indicators of a systematic development of metallicity in galaxies is doubtful, in quasars   at best non-existent and apparently already  refuted by empirical data.\footnote{\cite{Hasinger_ea:2002}, for a recent critical analysis see \cite{Yang_ea:age_crisis}. }
An interesting debate on the reliability of the  interpretation of the CMB anisotropy structure as indicator of primordial densitiy fluctuations has  started.\footnote{Cf. F. Steiner's contribution, this volume, and, among others, \cite{Ayaita_ea:CMB}.} 

It will be interesting to see, whether in the years to come developments of scalar field explanations in  cosmology harden and join in with new developments in high energy physics and observational cosmology; or whether they fall apart without contributing markedly to a better understanding of the challenging phenomena of elementary particle physics and cosmology. At  present there are good reasons to hope for  the first  outcome.


\bibliographystyle{apsr}
\bibliography{a_lit_hist,a_lit_mathsci}



\end{document}